\documentclass[onefignum,onetabnum]{siamart190516}

\usepackage{algorithm} 
\usepackage{algorithmic} 
\usepackage[algo2e,linesnumbered,boxed]{algorithm2e}

\usepackage{lipsum}
\usepackage{amsfonts}
\usepackage{graphicx}
\usepackage{epstopdf}
\usepackage{algorithmic}
\ifpdf
  \DeclareGraphicsExtensions{.eps,.pdf,.png,.jpg}
\else
  \DeclareGraphicsExtensions{.eps}
\fi

\newsiamremark{remark}{Remark}
\newsiamremark{hypothesis}{Hypothesis}
\crefname{hypothesis}{Hypothesis}{Hypotheses}
\newsiamthm{claim}{Claim}

\newsiamremark{example}{Example}
\headers{Companion Curve}{W. Wu and C. Chen}

\title{A Companion Curve Tracing Method for Rank-deficient Polynomial Systems}

\author{Wenyuan Wu\thanks{Chongqing Institute of Green and Intelligent Technology, Chinese Academy of Sciences, Chongqing, China
  (\email{wuwenyuan@cigit.ac.cn}, \email{chenchangbo@cigit.ac.cn},  \url{https://www.arcnl.org}).}
\and Changbo Chen\footnotemark[1] \thanks{Corresponding author.}}
\usepackage{amsopn}

\ifpdf
\hypersetup{
  pdftitle={An Example Article},
  pdfauthor={D. Doe, P. T. Frank, and J. E. Smith}
}
\fi

\DeclareMathAlphabet{\mathpzc}{OT1}{pzc}{m}{it}

\newcommand{\C}{\mathbb{C}}

\newcommand{\Jac}{\mathcal{J}}

\newcommand{\randpoint}{\ensuremath{\mathfrak{a}}}

\newcommand{\di}{\deg_{ind}}
\newcommand{\I}{\mathpzc{I}}
\newcommand{\R}{\ensuremath{\mathbb{R}}}

\newtheorem{prop}[theorem]{Proposition}
\newtheorem{cor}[theorem]{Corollary}
\newtheorem{define}[theorem]{Definition}

\def\proof{\textsc{Proof.} }
\def\foorp{\hfill$\square$}

\begin{document}

\maketitle


\begin{abstract}
  We propose a method for tracing implicit real algebraic curves defined by polynomials with rank-deficient Jacobians.
  For a given curve $f^{-1}(0)$, it first utilizes a regularization technique to compute at least one witness point per connected component of  the curve.
  We improve this step by establishing a sufficient condition for testing the emptiness of $f^{-1}(0)$.
  We also analyze the convergence rate and carry out an error analysis for refining the witness points.
  The witness points are obtained by computing the minimum distance of a random point to a smooth manifold embedding the curve while at the same time penalizing the residual
  of $f$ at the local minima.
  To trace the curve starting from these witness points,
  we prove that if one drags the random point  along a trajectory inside a tubular neighborhood of the embedded manifold of the curve,
  the projection of the trajectory on the manifold is unique and can be  computed by numerical continuation.
  We then show how to  choose such a trajectory to approximate the curve by computing eigenvectors of certain matrices.
  Effectiveness of the method is illustrated by examples.
\end{abstract}

\begin{keywords}
  rank-deficiency, real algebraic curve, curve tracing, numerical continuation, penalty function, {Tikhonov regularization}
\end{keywords}

\begin{AMS}
  65H10, 14Q30, 90C23
\end{AMS}

\section{Introduction}

Given an implicit real algebraic curve in $\R^n$ defined by a finite set of polynomials $f\subset\R[x_1,\ldots,x_n]$,
producing a polygonal chain approximation of it is a classical problem.
Existing numerical continuation methods~\cite{Allgower2003} for solving this problem often require that the Jacobian of $f$, denoted by $\Jac_f$, is of nullity one at all (or almost all) points of the real zero set of $f$,
that is $f^{-1}(0)$.
It remains a challenge to trace the curve defined by $f$ when $\Jac_f$ is rank-deficient, that is of nullity greater than one,
which is the dimension of the curve.
One typical example is when $f$ is a sum of squares of polynomials or more generally when $f$ is nonnegtive.

For a smooth curve in $\R^n$ defined by $f=\{f_1,\ldots,f_{n-1}\}$, where $\Jac_f$ is of full rank at any point of the curve,
the main technical challenge would be to identify all branches of $f^{-1}(0)$ and make sure that there is no jumping during curve tracing.
Identifying all branches of $f^{-1}(0)$ is the problem of computing witness points for every connected component of $f^{-1}(0)$~\cite{Rouillier2000,Hauenstein2012,WR13}.
Techniques for presenting or detecting curve jumping also exist~\cite{Blum1997,Beltran2013,Martin2013,Yu2014,WRF2017}.

For an almost smooth curve in $\R^n$ defined by $f=\{f_1,\ldots,f_{n-1}\}$, where $\Jac_f$ is of full rank at nearly all points of the curve,
one further difficulty is to trace across the singular points and get the correct topology around the singular points.
Many work exist for handling such problems~\cite{Hong1996,Daouda2008,Cheng2010,Gomes2014,DBLP:journals/jossac/ChenWF20,DBLP:journals/jossac/JinC20a}.

In this paper, we are interested in the case that $\Jac_f$ is rank-deficient at every point of $f^{-1}(0)$,
{ such as when $f$ consists of a polynomial in sum of squares.
  When $\Jac_f$ is rank-deficient, one of the main obstacles is that it is hard to find a tracing direction since
  the dimension of the nullspace of $\Jac_f$ is at least two.
}
Our starting point for solving this problem is a penalty function based method for computing witness points for every connected component of $f^{-1}(0)$~\cite{Wu2017}.
We then extend it to a companion curve tracing method.
The main idea of the penalty function is to embed the curve in a high-dimensional smooth manifold defined by a system $g$
with full rank Jacobians and compute the minimum distance of a random point to the manifold while at the same time penalizing the residual
of $f$ at the local minima.
{ 
To render the curve, one natural idea is to move the random point (as a guiding point)  along a trajectory and hopefully
the corresponding minima will be dragged continuously.
To implement this idea, there are two main challenges to overcome.
  One is the potential occurence of discontinuity, which indeed may happen
  as illustrated in Section~\ref{sec:intro-ex} and proved in Section~\ref{sec:accuracy}.
Another is to make sure the minima do move along the curve as one drags the guiding point along the trajectory.
We show that both challenges can be overcome and the idea of moving the guiding point is indeed feasible if the  point is moved along the directions
defined by eigenvectors of certain matrices
and inside a tubular neighborhood of the embedded manifold of the curve.
}
The trajectory formed by moving the random point is called a {\bf companion curve} of $f^{-1}(0)$.

Interestingly,  the  penalty function method, although initially proposed in a completely different context, 
turns out to be closely related to Tikhonov regularization for solving rank-deficient nonlinear least squares problems~\cite{Tikhonov95,Engl1996,Eriksson2005},
which is explained in the preliminary section.

The paper is structured as follows.
In Section~\ref{sec:pre}, we recall how to compute at least one witness point
for every connected component of an arbitrary real algebraic variety,
whose defining system allows to be rank-deficient, via the so-called penalty function method~\cite{Wu2017}.
In Section~\ref{sec:intro-ex}, 
we illustrate by a simple example the challenge for tracing real algebraic curves defined by rank-deficient systems,
with the initial witness points provided by the penalty function method, as well as
the main idea of our companion curve method for handling this problem.
One weakness of the penalty function method for computing witness points is that
it always return a non-empty set of points even the given real variety is empty.
In Section~\ref{s:emptiness}, we provide a sufficient criterion for testing emptiness of
a variety by the penalty function method.
To successfully tracing a real algebraic curve, it is important to make the approximate witness points
close enough to the curve.
In Section~\ref{sec:accuracy}, we propose a homotopy method for improving the precision of witness points.
With all these preparations and another tool from differential geometry,
namely Tubular Neighborhood Theorem, we propose a companion curve method
for curve tracing in  Section~\ref{sec:pathtrack}.
{ The effectiveness of the method is illustrated by several examples in Section~\ref{sec:ex}.}
Finally, in Section~\ref{sec:con}, we draw the conclusion
and propose several ways to improve the current method.

\section{Preliminaries}
\label{sec:pre}
In this section, we recall some preliminary results
that were introduced in~\cite{Wu2017} for computing witness points of
rank-deficient polynomial systems.

Let $x=(x_1,\ldots,x_n)$ and let $f = \{f_1,...,f_k\}\subset \R[x]$.
Let $V_{\R}(f)$ be the zero set of $f$ in $\R^n$.
Let $\randpoint=(\randpoint_1,...,\randpoint_n)\notin V_{\R}(f)$ be a point in $x$-space and consider the minimal distance
from $V_{\R}(f)$ to this point:
\begin{eqnarray}\label{eq:opt2}
 \min \;  \sum_{i=1}^n (x_i-\randpoint_i)^2 \\
 s.t. \hspace{1cm} f(x) = 0   \nonumber.
\end{eqnarray}
Clearly every semi-algebraically connected component of $V_{\R}(f)$
has at least one point attaining the local minimum.
These points make the following matrix lose full rank:
$$
A := \left(
\begin{array}{ccc}
              \partial f_1/\partial x_1 & \cdots & \partial f_1/\partial x_n \\
              \vdots & \ddots & \vdots \\
              \partial f_k/\partial x_1 & \cdots & \partial f_k/\partial x_n \\
              x_1-\randpoint_1 & \cdots & x_n-\randpoint_n\\
            
\end{array}
\right).
$$
Note that the first $k$ rows of $A$ are exactly the the Jacobian matrix of $f$ w.r.t. $x$, denoted by $\Jac_f$.
If $\Jac_f$ is rank-deficient at every point of $V_{\R}(f)$, the matrix $A$ automatically loses
full rank and thus does not provide any extra helpful information on computing the semi-algebraically connected components of $V_{\R}(f)$.

To overcome this algebraic rank deficiency problem, the paper~\cite{Wu2017} introduces a penalty function based approach.
Instead of solving the optimization problem~(\ref{eq:opt2}), one considers the following unconstrained optimization problem:
\begin{eqnarray}\label{eq:optu}
 \min \; \mu &= ( \beta\cdot (f_1^2+\cdots+f_k^2) +  \sum_{i=1}^n (x_i-\randpoint_i)^2) /2.
\end{eqnarray}
Note that as $\beta$ approaches infinity, $f_i$, $i=1,\ldots, k$ are forced to be zero.
Intuitively, this provides an approximate solution to the problem~(\ref{eq:opt2}) for large enough $\beta$.
Indeed, the paper~\cite{Wu2017} proves the following result justifying such intuitions.
\begin{prop}[Corollary $1$ in~\cite{Wu2017}]
  \label{lem:inf}
Let $p$ be a local minimum of (\ref{eq:opt2}). There exists a local minimum $p'$ of (\ref{eq:optu}) for sufficiently large $\beta$, such that $\|p-p'\|$ can be arbitrarily small.
\end{prop}
Moreover, one can get a rough estimation of the distance between $p$ and $p'$ via the notion of {\em degree index}.

\begin{define}[Definition $3$ in~\cite{Wu2017}]
  \label{def:di}
  For a given $v\neq 0 \in \R^n$, let $f_v := f(x=vt)$.
  Denote by $\deg_{\min}(f_v)$ the trailing degree of $f_v$.
  We call $\deg_{ind}(f)  = \max_{v} \deg_{\min}(f_v)$ the {\em degree index} of $f$.
  Given a point $p\in V_{\R}(f)$, the degree index of $f$ at $p$ is defined as $\deg_{ind}(f(x+p))$.
 \end{define}

\begin{theorem}[Theorem $5$ in~\cite{Wu2017}]
  \label{thm:est2}
  For a random point $\randpoint \in \R^n$ and a sufficiently large $\beta$, suppose that $p\in V_{\R}(f)$ attains the local minimal distance to $\randpoint$.
  Then there is a solution $p'$ of Equation (\ref{eq:newsys3}) such that $\|p'-p\| = O( \sqrt[2\I-1]{1/\beta} \, )$, where $\I = \max\{\di(f_i(x+p)),i=1,...,k\}$.
\end{theorem}

The local minima of (\ref{eq:optu}) are exactly the points that vanish the gradient of $\mu$, that is satisfying the following equation:
\begin{equation}\label{eq:newsys3}
\left(
  \begin{array}{c}
    x_1 \\
    \vdots \\
    x_n \\
  \end{array}
\right) + \beta \cdot \Jac^t \cdot \left(
                         \begin{array}{c}
                           f_1 \\
                           \vdots \\
                           f_k \\
                         \end{array}
                       \right) =  \left(
                         \begin{array}{c}
                           \randpoint_1 \\
                           \vdots \\
                           \randpoint_n \\
                         \end{array}
                       \right).\\
\end{equation}
where the $n\times k$ matrix $\Jac^t$ is the transpose of the Jacobian of $f$.

The left hand side of  Equation (\ref{eq:newsys3}) defines a smooth mapping $M:\R^{n}\rightarrow\R^n$.

\begin{lemma}[Lemma $2$ in~\cite{Wu2017}]
  \label{lem:regular}
  For almost all points $ \randpoint = (\randpoint_1,...,\randpoint_n)\notin V_{\mathbb{R}}(f)$, $M^{-1}(\randpoint)$
  is a nonempty finite set and every point of  $M^{-1}(\randpoint)$ is a regular point of $M$.
\end{lemma}

Proposition~\ref{lem:inf} and Lemma~\ref{lem:regular} together show that for almost all points $\randpoint$ of $\R^n$, $M^{-1}(\randpoint)$
  contains points meeting every semi-algebraically connected component of $V_{\R}(f)$.
  We warn that $M^{-1}(\randpoint)$ may contain extra points not belonging to $V_{\R}(f)$.
  In particular, even $V_{\R}(f)$ is empty, $M^{-1}(\randpoint)$ is always nonempty.
  Numerically testing if a real variety is empty in general is a difficult problem.
  We provide a partial answer to this problem in Section~\ref{s:emptiness}.

  Sometimes, it is useful to use the following two equivalent formulations to (\ref{eq:optu}).

  Let $z=(z_1,..,z_k)$ be $k$ slack variables and $g=\{f_1+z_1, f_2+z_2, ..., f_k+z_k\}$.
Note that we have $g\subset \R[x, z]$ and $V_{\R}(g)\subseteq \R^{n+k}$.
  \begin{eqnarray}\label{eq:opt}
 \min \; \mu &= ( \beta\cdot (z_1^2+\cdots+z_k^2) +  \sum_{i=1}^n (x_i-\randpoint_i)^2) /2 \\
 &s.t. \hspace{1cm} g(x,z)=0  \nonumber.
\end{eqnarray}

Let $z_i=w_i/\sqrt{\beta}$, $i=1,\ldots,k$,
and substitute them into (\ref{eq:opt}).
Let $h=\{f_1+w_1/\sqrt{\beta},\ldots,f_k+w_k/\sqrt{\beta}\}$.

\begin{eqnarray}\label{eq:optd}
 \min \; & ( w_1^2+\cdots+w_k^2 +  \sum_{i=1}^n (x_i-\randpoint_i)^2) /2 \\
 &s.t. \hspace{1cm} h=0 \nonumber.
\end{eqnarray}

One nice thing about (\ref{eq:opt}) and (\ref{eq:optd}) is that
the real varieties defined by $g$ and $h$ are smooth submanifolds of $\R^{n+k}$.

\begin{remark}
  \label{remark:regularization}
If we replace $\beta=1/t$, we obtain another equivalent formulation:
\begin{eqnarray}\label{eq:optr}
 \min \; \mu &= ((f_1^2+\cdots+f_k^2) +  t\sum_{i=1}^n (x_i-\randpoint_i)^2) /2.
\end{eqnarray}
Such a formulation is exactly the Tikhonov regularization for nonlinear least squares problems~\cite{Tikhonov95,Engl1996,Eriksson2005},
where $t\sum_{i=1}^n (x_i-\randpoint_i)^2) /2$ is the regularization term for the minimization problem: $\min \|f\|^2$.
{ Here, one must be cautious about the choice of $\randpoint$, since Lemma~\ref{lem:regular} does not exclude
  the possibility that  the solution of problem (\ref{eq:optr}) is not regular for certain $\randpoint$, which indeed poses a challenge for curve tracing as explained in next section.}
\end{remark}

\section{An introductory example}
\label{sec:intro-ex}
In this section, we illustrate by an example the challenge of a pure numerical method
for tracing algebraic curves defined by rank-deficient polynomials as well as
the main idea of our companion curve method.



\begin{figure}
    \centering
 \begin{minipage}[c]{0.55\linewidth}
\centering
\includegraphics[width=\textwidth]{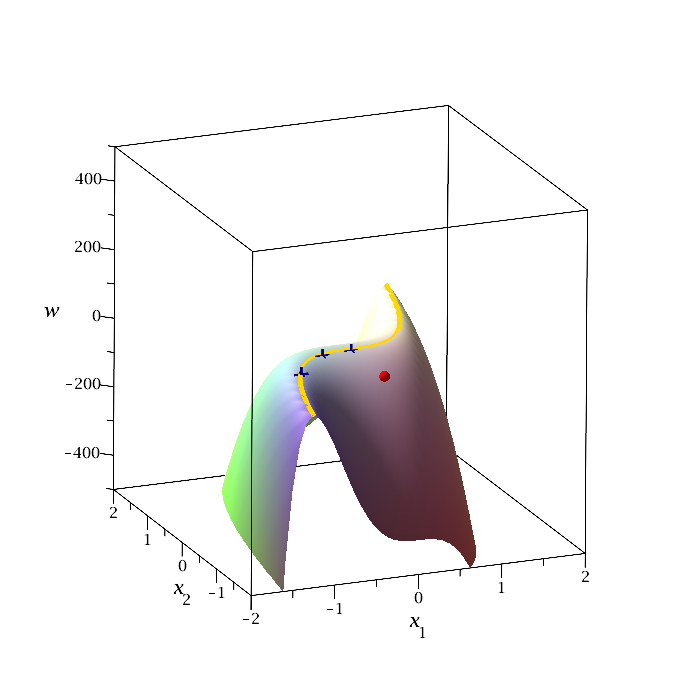}
\end{minipage}
\begin{minipage}[c]{0.43\linewidth}
\centering
\includegraphics[width=\textwidth]{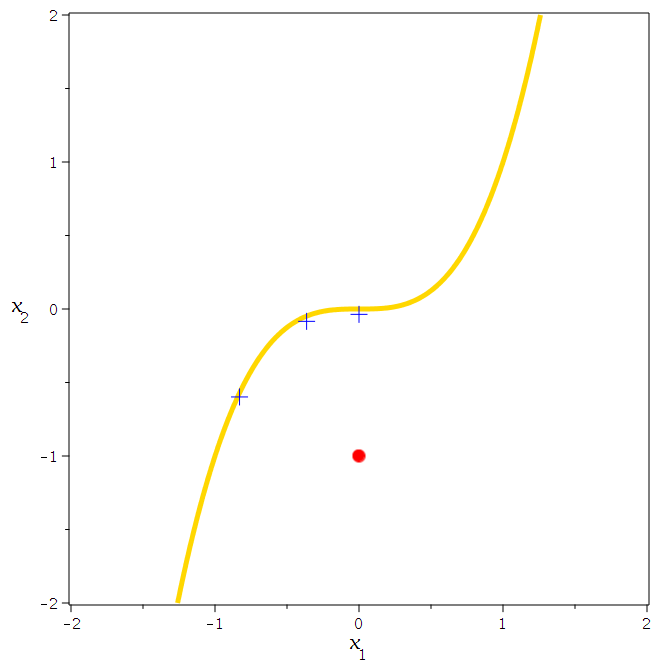}
\end{minipage}
\caption{{(Color online) Left: the random point $\randpoint=(0, -1, 0)$ (in red $\bullet$), the three local minima (in blue $+$),
    the smooth surface defined by $f+w/100=0$ and its intersection with $w=0$ (gold curve).
    Right: the random point $\randpoint=(0, -1)$ (in red $\bullet$), the three approximate witness points (in blue $+$),
    and the curve $f=0$ (gold curve). }}
\label{fig:cubic-wp}
\end{figure}

Let $f:={x_{{1}}}^{6}-2\,{x_{{1}}}^{3}x_{{2}}+{x_{{2}}}^{2}=(x_1^3-x_2)^2$.
Recall from Section~\ref{sec:pre} that, one can obtain approximate witness points of $V_{\R}(f)$ 
by solving Equation (\ref{eq:newsys3}).
Choosing $\beta=10^4$ and $\randpoint=(0, -1)$, the equation becomes $\{60000\,{x_{{1}}}^{11}-180000\,{x_{{1}}}^{8}x_{{2}}+180000\,{x_{{1}}}^{5}{x_{{2}}}^{2}-60000\,{x_{{1}}}^{2}{x_{{2}}}^{3}+x_{{1}},
-20000\,{x_{{1}}}^{9}+60000\,{x_{{1}}}^{6}x_{{2}}-60000\,{x_{{1}}}^{3}{x_{{2}}}^{2}+20000\,{x_{{2}}}^{3}+x_{{2}}+1\}$.
By the homotopy continuation method, say the one implemented in Hom4PS-2.0~\cite{Lee2008}, one obtains three approximate  witness points: $(-0.3639,-0.0840)$, $(-0.8296,-0.5982)$, $(0,-0.0364)$.
Geometrically, these  points are actually the projection onto $(x_1,x_2)$-space of the following three local minima of the optimization problem (\ref{eq:optd}):
$$(-0.3639,-0.0840, -0.1280), (-0.8296,-0.5982, -0.0739), (0,-0.0364, -0.1324).$$
These three points attain the local minimum distance from the point $(0, -1, 0)$
to the manifold $f+w/100=0$, as illustrated by Figure~\ref{fig:cubic-wp}.

\subsection{The challenge}
\begin{figure}
    \centering
 \begin{minipage}[c]{0.48\linewidth}
\centering
\includegraphics[width=\textwidth]{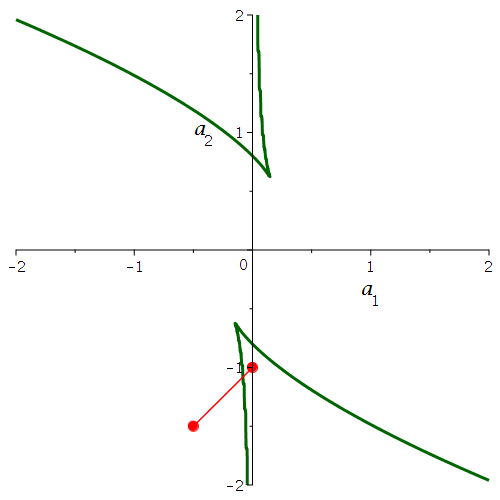}
\end{minipage}
\begin{minipage}[c]{0.48\linewidth}
\centering
\includegraphics[width=\textwidth]{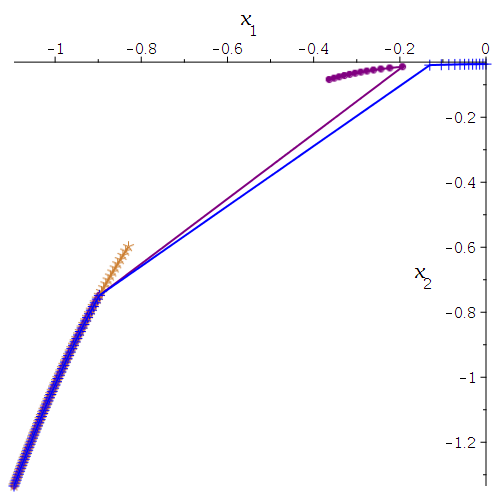}
\end{minipage}
\caption{{(Color online) Left: the discriminant locus of $H(x, \randpoint)$ (two ``V'' curves in dark green) and the segment $\overline{p_0p_1}$ (in red).
    Right: the traced curve starting with the initial point $(-0.8296,-0.5982)$ (gold $\star$), the traced curve starting with the initial point $(0, -0.0364)$ (blue $+$),
    and the traced curve starting with the initial point $(-0.3639,-0.0840)$ (purple $\bullet$).}}
\label{fig:discrim}
\end{figure}

Next we would like to generate more points of $V_{\R}(f)$ by curve tracing with the three initial points.
One natural idea is to consider the following linear homotopy:
\begin{equation}\label{eq:tau}
H_{\tau}(x,\tau)= x+\beta_0 \Jac^t\cdot f - (\tau p_1 + (1-\tau) p_0) \equiv 0,
\end{equation}
where $\beta_0=10000$ and one moves $\randpoint$ from $p_0=(0, -1)$ to $p_1=(-0.5,-1.5)$ in a line.
As long as the number of solutions of $H(x, \randpoint)=x+\beta_0 \Jac^t\cdot f -\randpoint$ in $x$ remains
unchanged and the graphs of the solutions of $H(x, \randpoint)$ as functions of $\randpoint$ remain disjoint and smooth,
one should encounter no much difficulty during curve tracing.
However, as illustrated by the left subfigure of Fig.~\ref{fig:discrim}, the path of $\randpoint$ crosses
the discriminant locus~\cite{LazardRouillier2007} of $H(x, \randpoint)$,
which can be obtained by a Gr\"obner basis computation~\cite{LazardRouillier2007} by Lemma~\ref{lem:nonzeropoly} in Section~\ref{sec:accuracy}.
As illustrated by the right subfigure of Fig.~\ref{fig:discrim},
when $\randpoint$ approaches the discriminant locus, the solution curve starting with $(0, -0.0364)$ (blue $+$) and the solution curve starting with $(-0.3639,-0.0840)$ (purple $\bullet$)
get close to each other gradually and finally are no longer solutions of $H(x, \randpoint)$.
Interestingly, if we apply Newton iteration now starting with these non-solution points,
they finally converge to the third solution curve traced starting with the  local minimum $(-0.8296,-0.5982)$ (gold $\star$).
Therefore, curve jumping happened when the parameter $\randpoint$ crossed the discriminant locus.

\subsection{The companion curve solution}

Note that in Fig.~\ref{fig:cubic-wp}, there are three points on the smooth surface attaining the local minimum distance to the given random point $(0, -1, 0)$.
Intuitively, there should only be one local minimum point if the given random point is close enough to the surface and the surface is connected.
This is indeed true, which will be proved in Section~\ref{sec:pathtrack}, thanks to the Tubular Neighborhood Theorem in differential geometry.

\begin{figure}
    \centering
 \begin{minipage}[c]{0.48\linewidth}
\centering
\includegraphics[width=\textwidth]{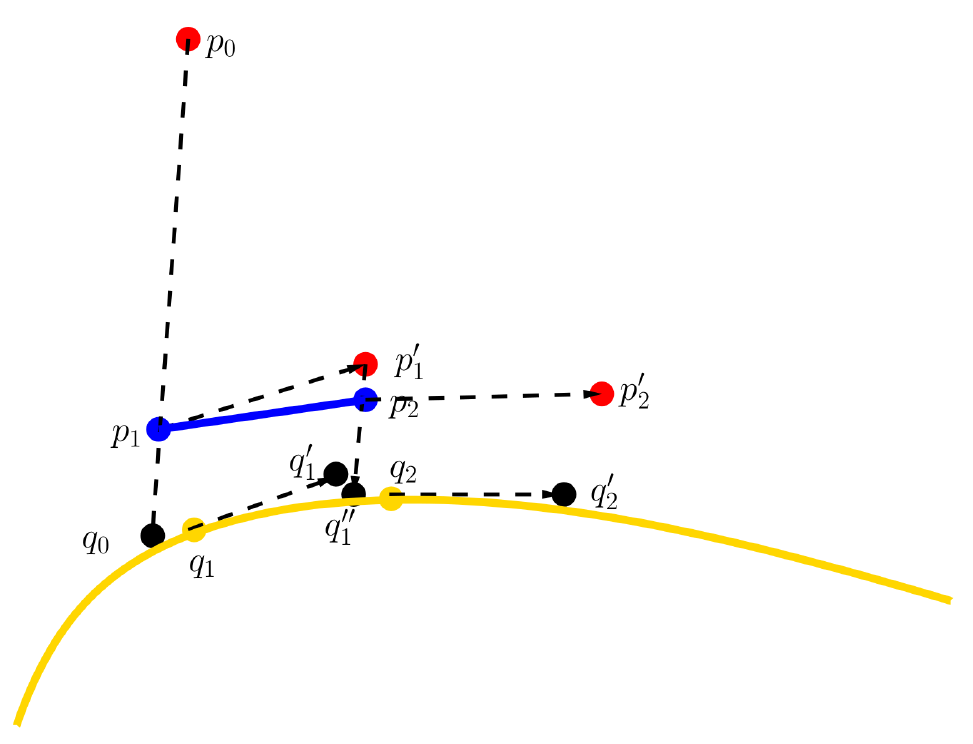}
\end{minipage}
\begin{minipage}[c]{0.48\linewidth}
\centering
\includegraphics[width=\textwidth]{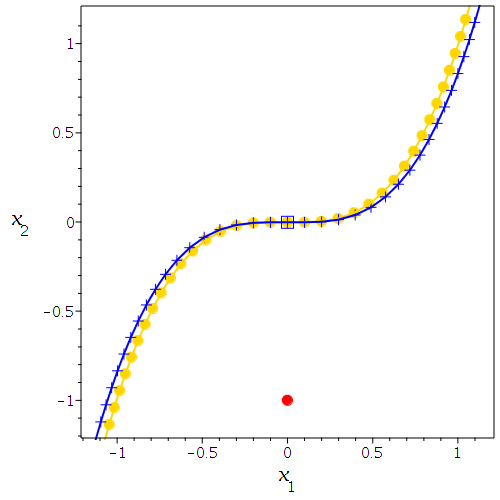}
\end{minipage}
\caption{(Color online) Left: illustrating the idea of companion curve tracing method.
    Right: an approximation of $V_{\R}(f)$ (polygonal chain in gold $\bullet$) and its companion curve (polygonal chain in blue $+$).
  The initial random point $(0, -1)$ (in red $\bullet$) at the bottom was moved to a point (in blue $\Box$) close to $V_{\R}(f)$.}
\label{fig:cubic-comp}
\end{figure}

Now suppose that $\randpoint=p_0$ is the initial random point and $(x,w)=(q_0,r_0)$ is an approximate local minimum of the optimization problem (\ref{eq:optd}).
Suppose that $\|r_0\|$ is small enough such that $x=q_0$ can be seen as an approximate witness point of $V_{\R}(f)$.
We first move the point $\randpoint$ to another point $p_1$ on the segment $\overline{p_0q_0}$ such that $(p_1, 0)$ is close to the surface and hope that there is only one point $(q_1, r_1)$
on the surface with local minimum distance to $(p_1, 0)$ nearby. We then pick a well chosen direction, say $\vec{v}$,
{ by some eigenvectors computation and move $\randpoint$ from $p_1$ to $p_1'$.
  Accordingly, $x$ is moved in the same direction from $q_1$ to $q_1'$ (may be further refined to $q_1''$ by Newton iteration if necessary).
}
It is possible that $(p_1',0)$ is now outside the tubular neighborhood of the surface and one may repeat the previous step to
drag $p_1'$ to $p_2$ and produce $q_2$. And then move $q_2$ to $q_2'$, $p_2$ to $p_2'$, so on so forth.
The polygonal chain formed by the sequence of points $p_1, p_2, \ldots$ is called the companion curve of $V_{\R}(f)$,
as illustrated in  Fig.~\ref{fig:cubic-comp}.

\section{Emptiness of real variety}
\label{s:emptiness}
In this section, we propose a criterion for testing the emptiness of a given real variety.

Lemma \ref{lem:regular} implies that all the solutions of Equation (\ref{eq:newsys3})
can be obtained by applying homotopy continuation methods.
Among these solutions, we look for solutions with small residuals i.e. $\|z\| \ll 1$.
It is possible that such points do not exist, which then provides strong evidence that $V_{\R}(f)$ is empty.
Intuitively, this is because if $V_{\R}(f)$ is not empty, increasing the penalty factor $\beta$
will force $\|z\|$ close to zero. Thus, the minimal value of $\mu$ will be slightly larger than the distance from
$\randpoint$ to $V_{\R}(f)$.

To study the relationship between $\mu_{\min}$ and the emptiness of $V_{\R}(f)$,  we homogenize the system $f$ by adding a variable $x_0$ satisfying a new equation $\bar{f}_{k+1} = \sum_{i=0}^n x_i^2 - 1 = 0$ to obtain a \textit{homogenized system} $\bar{f}$ except for the new inhomogeneous equation. The corresponding unconstrained optimization problem is
\begin{eqnarray}\label{eq:optmu2}
 \min \; \bar{\mu} &= ( \beta\cdot ( \bar{f}_1^2+\cdots+\bar{f}_k^2 + \bar{f}_{k+1}^2 ) +  \sum_{i=0}^n (x_i-\randpoint_i)^2) /2.
\end{eqnarray}
where $\randpoint$ is  chosen randomly in the unit ball $B(0;1)$ of $\R^{n+1}$.

\begin{prop}\label{prop:empty}
Let $\bar{\mu}_{\min}$ be the global minimal value of $\bar{\mu}$ in the optimization problem (\ref{eq:optmu2}).
If $\bar{\mu}_{\min} > 2$, then $V_{\R}(f)= \emptyset$.
\end{prop}
\proof
We prove it by contradiction. Suppose $V_{\R}(f) \neq \emptyset$. Then $V_{\R}(\bar{f}) \neq \emptyset$ and let $p \in V_{\R}(\bar{f})$. We know that $p,\randpoint \in B(0;1)$ which implies $\|p-\randpoint\| \leq 2$. Thus, $\bar{\mu}_{\min} \leq \beta\cdot 0+ \|p-\randpoint\|^2/2 \leq 2$. It contradicts the assumption $\bar{\mu}_{\min} > 2$.
\foorp\\

\begin{example}
 Let $f= {x}^{4}+{y}^{4}-3\,xy+3/2$.  We can verify that $$f={\frac {159}{784}}\,{x}^{4}+{\frac {15}{64}}\,{y}^{4}+{\frac {3}{50}}+
 \left( \frac{6}{5}-\frac{5}{4}\,xy \right) ^{2}+ \left( {\frac {25}{28}}\,{x}^{2}-{
\frac {7}{8}}\,{y}^{2} \right) ^{2} >0.
$$
Homogenizing $f$ yields $\bar{f}=\{{x}^{4}+{y}^{4}-3\,{h}^{2}xy+3/2\,{h}^{4},{h}^{2}+{x}^{2}+{y}^{2}-1\}$.
Choose $\randpoint=(0.2,0.5,0.3)$ and $\beta=10000$ and solve the corresponding system of Equation (\ref{eq:newsys3}) by \texttt{Hom4Ps2}. It gives $23$ real roots among $111$ complex ones. The minimal value of $\bar{\mu}$ is $28.6$ at $(x=0.565,y=0.565,h=0.596)$.
By Proposition \ref{prop:empty}, it indicates that $V_{\R}(f)= \emptyset$.
\end{example}

To apply this lemma we may choose sufficiently large $\beta$ to make $\bar{\mu}_{\min} > 4$ possible. However, for some positive polynomials it will never happen no matter how large  $\beta$ is.

\begin{example}
Consider $f= (xy-1)^2+y^2$ and it is positive but arbitrarily close to zero.
Fixing $\randpoint=(0,0,0)$, by Equation (\ref{eq:optmu2}), we have
$$\bar{\mu} = \beta \left( {h}^{4}-2\,{h}^{2}xy+{h}^{2}{y}^{2}+{x}^{2}{y}^{2}
 \right) ^{2}+ \beta \left( {h}^{2}+{x}^{2}+{y}^{2}-1 \right) ^{2}+{h
}^{2}+{x}^{2}+{y}^{2}.
$$
Numerical computation shows that $\bar{\mu}_{\min}=0.999975$ at  the point $(x = 2.9\times 10^{-8}, y = 0.9999, h = 8.5\times 10^{-9})$  when $\beta=10^4$.
Actually, $\bar{\mu}_{\min}$ must be no greater than $1$ since $\bar{\mu}(0,1,0) = 1$ for any $\beta$. It means that we cannot tell if
the real variety of $f$ is empty or not by Proposition \ref{prop:empty}.
\end{example}

Proposition \ref{prop:empty} only gives a sufficient condition for $V_{\R}(f)= \emptyset$. If $\bar{\mu}_{\min} < 2$,
{ deciding the emptiness is an open question.}
In the rest of this paper, we always assume that $V_{\R}(f)\neq\emptyset$.



\section{Refinement}
\label{sec:accuracy}


By Proposition \ref{lem:inf}, theoretically we can use Equation (\ref{eq:newsys3}) to update the approximate root $x'$ by increasing $\beta$. Suppose we have all the real solutions of Equation (\ref{eq:newsys3}) for $\beta=\beta_0$ denoted by $R_{\beta_0}$.
When we increase $\beta$ from $\beta_0$ to $\beta_1$, there are two ways to obtain $R_{\beta_1}$. One way is to solve Equation (\ref{eq:newsys3}) in the complex field and then keep the real solutions. Alternatively, we may trace the real curves
of the following homotopy starting from
all points of $R_{\beta_0}$ by moving $\beta$ from $\beta_0$ to $\beta_1$ continuously.
\begin{equation}\label{eq:homotopy1}
 H(x,\beta)= (x-\randpoint) + \beta\cdot \Jac^t \cdot f\equiv 0.
\end{equation}

But we have to be aware of singular Jacobian for a successful tracing.

\begin{lemma}\label{lem:nonzeropoly}
For any polynomial system $f\subset  \R[x]$, let $F = \{{x}-\randpoint + \beta\cdot \Jac^t \cdot {f}\}$ and $G= \{F, \det(\frac{\partial F}{\partial x})  \} \subset \R[x,\randpoint,\beta]$. Then there is a nonzero polynomial $\phi(\randpoint,\beta) \in \langle G \rangle$.
\end{lemma}
\proof
When $\langle G \rangle = \langle 1 \rangle$, it is trivial.\\
Otherwise, let $Gb$ be a Gr\"obner basis of $G$ with respect to lex order $x\succ \randpoint \succ \beta$.
Then $Gb' = Gb \cap \R[\randpoint,\beta]$ is a Gr\"obner basis of the elimination ideal.
If $Gb' = 0$, then for a generic $(\randpoint,\beta) \in \C^{n+1}$, it can be extended to a solution $z$ of $G=0$.
But by Sard's lemma for varieties (Chap. 3 in \cite{Mumford1981}), this point is a regular value of the smooth mapping with nonsingular Jacobian. Thus, $z$ does not satisfies $\det(\frac{\partial F}{\partial x}) =0$ which contradicts the assumption $z$ is a solution of $G=0$.

Therefore, $Gb'$ must contain nonzero polynomials.
\foorp\\

This lemma shows that even $\randpoint$  is chosen randomly, it is still possible to encounter singular Jacobian when $\beta$ increases continuously from $\beta_0$ to infinity.

\begin{example}
  Recall the example $f:=(x_1^3-x_2)^2$ in Section~\ref{sec:intro-ex}.
  If we choose $\randpoint=(0,-1)$, Equation (\ref{eq:newsys3}) becomes
  $$
  \left\{
  \begin{array}{rcl}
  6\,\beta\,{x}^{11}-18\,\beta\,{x}^{8}y+18\,\beta\,{x}^{5}{y}^{2}-6\,
  \beta\,{x}^{2}{y}^{3}-{\it a1}+x&=&0\\
  -2\,\beta\,{x}^{9}+6\,\beta\,{x}^{6}y-6\,\beta\,{x}^{3}{y}^{2}+2\,\beta\,{y}^{3}-{\it a2}+y\}&=&0
  \end{array}
  \right..
  $$
  Applying Lemma~\ref{lem:nonzeropoly}, we obtain
  $
  \phi(\randpoint,\beta)=2082930190011\,{\beta}^{5}-47659837219452\,{\beta}^{4}-19398995284788\,{\beta}^{3}-3070008000000\,{\beta}^{2}-207900000000\,\beta-5000000000,
  $
  which has only one positive real root $\beta\approx23.28$.
  So if we choose $\beta_0\leq 23$, singular Jacobian will be encountered when increasing $\beta$ from $\beta_0$ to infinity. 
\end{example}

Next we will study such probability.
After we fix the value of $\randpoint$, $Gb'$ is a set of univariate polynomials and it is generated by a single polynomial denoted by $g \in \R[\beta]$. By rescaling, the all coefficients of $g$ are in $[-1,1]$.

\begin{lemma}\label{lem:prob}
Let $g= \sum_{i=0}^n c_i \beta^i \in \R[\beta]$ with degree $n$. Suppose the coefficients $\{c_0,...,c_n\}$ are i.i.d. random variables with uniform distribution in $[-1,1]$. Let $k$ be the largest absolute value of all the roots of $g$. Then the probability  $\Pr(k<N+1) > 1- 1/N$.
\end{lemma}
\proof
By the root bound in Chap. 8 of \cite{Mishra93}, we have  $$k < \frac{|c_n|+\max(|c_{n-1}|,...,|c_0|)}{|c_n|}.$$
Then clearly, it implies $\Pr(k<N+1) > \Pr(\max(|c_{n-1}|,...,|c_0|) \leq |c_n|N )$.

Let $c= \max(|c_{n-1}|,...,|c_0|)$ which is a random variable in $[0,1]$.
Moreover $\Pr(c\leq t) = t^n$ if $t<1$,  and $\Pr(c\leq t) = 1$ if $t\geq 1$.
Hence
\begin{eqnarray*}
  \Pr(c \leq |c_n|N ) &=& \Pr(|c_n|\geq 1/N) + \Pr(|c_n|< 1/N  \wedge c \leq |c_n|N )\\
                    &=& 1-1/N + \int_{0}^{1/N}(|c_n|N)^n~ d|c_n|  \\
                    &=& 1-1/N +  { 1/(n+1)/N} > 1-1/N.
\end{eqnarray*}
It gives $\Pr(k<N+1) >1- 1/N$.
\foorp \\

This lemma indicates that the probability of singularity occurring during the homotopy path tracking (\ref{eq:homotopy1}) is less than $1/(\beta_0 -1)$. Hence we usually choose a large $\beta$ e.g. $\beta=10^4$. \\

Proposition~\ref{lem:inf} shows the existence of an approximate point $p'$ for any critical point $p$ by  solving Equation (\ref{eq:newsys3}).
With the distribution assumption in Lemma \ref{lem:prob}, we can also show the uniqueness of $p'$ with a quite high probability.

\begin{cor}\label{cor:est3}
For a random point $\randpoint \in \R^n$ and a sufficiently large $\beta$, suppose $p\in V_{\R}(f)$ attains the local minimal distance to $\randpoint$.  Then there is a unique solution $p'$ of Equation (\ref{eq:newsys3}) such that $\|p'-p\| = O( \sqrt[2\I-1]{1/\beta} \, )$, where $\I = \max\{\di(f_i(x+p)),i=1,...,k\}$ with probability at least $1-1/(\beta-1)$.
\end{cor}
\proof
Compared with Theorem \ref{thm:est2}, we only need to show the uniqueness of $p'$.
Suppose there is another point $p''$ also close to $p$ such that $\|p''-p\| = O( \sqrt[2\I-1]{1/\beta} \, )$.

By Lemma \ref{lem:prob}, when $\beta\rightarrow \infty$,  with probability at least $1-1/(\beta-1)$,  both $p'$ and $p''$ approach $p$
as $\beta$ moves to $\infty$ continuously. And it means that for a generic $\randpoint$, $\phi(\randpoint,\beta) \rightarrow 0$  when $\beta \rightarrow \infty$, where $\phi$ is a nonzero polynomial given in
Lemma \ref{lem:nonzeropoly}. This indicates that $g(\randpoint)= 0$ where $g$ is the leading coefficient of $\phi$ with respect to $\beta$.
But $\randpoint$ is generic,  $g(\randpoint)\neq 0$. It leads to a contradiction. Therefore, $p'$ is unique for sufficiently large $\beta$.
\foorp \\

Replacing $\beta$ with $1/t$ in (\ref{eq:homotopy1}), we obtain a homotopy
\begin{equation}
  \label{eq:homotopy2}
   H(x, t) := t(x-\randpoint) + \Jac^t \cdot f  = 0.
\end{equation}
Let $(t_0, x_0)$ be an initial point satisfying $H(x_0, t_0)=0$.
As $t$ approaches zero, with the assumption in Lemma~\ref{lem:prob}, the homotopy path $x(t)$ will approach $V_{\R}(f)$ in high probability (about $(1-t)$).
As a direct consequence of Theorem \ref{thm:est2}, we have the following error estimation for the refinement process.
\begin{corollary}
 Let $d$ be the degree of $f$.
 If we reduce the value of $t$ by a half at each step and assume that the step sizes $\frac{t}{2^s}$, $s=1,2,\ldots$
 are small enough to avoid curve jumping, then after $s$ steps of path tracking,
 the error of root is reduced to $O(\tau^s\, \delta)$, where $\tau = 2^{-1/(2d-1)}$, and $\delta$ is the initial error $\|p'-p\|$.
\end{corollary}
\begin{remark}
This result was stated as  Corollary $6$ in~\cite{Wu2017} without probability discussion.   
\end{remark}

One could also estimate the backward error.
We can consider the embedding system $g:=\{f_1+z_1, \ldots, f_k+z_k\}$ as a perturbed system $f'$ of the input $f$ when fixing the values of $z_i$.
More precisely, a critical point $p'$ which is a solution of Equation (\ref{eq:newsys3}) is an exact solution of $f' = \{f_1+z_1,...,f_k+z_k\}$
where $z_i= -f_i(p')$. By Theorem~\ref{thm:est2}, we have the following result.

\begin{cor}
Let $p'$ be a local minimum of $\mu$ in~(\ref{eq:opt}) with a sufficiently large $\beta$. Then the residual error
satisfies
\begin{equation}\label{eq:res}
{\sum z_i^2 }= O(1/\beta^{\frac{2\I+1}{2\I-1}}).
\end{equation}
\end{cor}
\proof
Since $\mu(p')$ is a local minimal value, then the perturbations $z_i$ are very small because
\begin{equation*}\label{eq:z}
 \mu(p') \leq \mu(p) \Rightarrow \beta (\sum z_i^2) + \|p'-\randpoint\|^2 \leq \|p-\randpoint\|^2   \Rightarrow \beta (\sum z_i^2) \leq \|p'-p\|^2
\end{equation*}
where $p$ is a solution of Problem (\ref{eq:opt2}) and the forward error $\|p'-p\| = O( \sqrt[2\I-1]{1/\beta} \, ) $ when $\beta$ is sufficiently large by Theorem~\ref{thm:est2}.  It implies that
the backward error ${\sum z_i^2 }= O(1/\beta^{\frac{2\I+1}{2\I-1}})$.
\foorp

\section{Tracing algebraic curves defined by rank-deficient systems}
\label{sec:pathtrack}

Previously we have presented a method to approximate and refine the real witness points of a general system
combining critical point and penalty function techniques.
Further applications may require to produce more points on the real variety, especially when the variety is an algebraic curve.
It is  quite challenging for rank-deficient systems, since the Jacobian is always singular along the whole curve.
A straightforward method for generating more points by moving the random point $\randpoint$
to generate more critical points may fail as illustrated in Section~\ref{sec:intro-ex}.

Thanks to one of the fundamental theorems
in differential geometry: A smooth manifold is contained in an
open tubular neighborhood in which every point can be uniquely
projected onto the manifold following a normal line, we can track the critical point along the curve if $\randpoint$ is always contained
in this tubular neighborhood.

\subsection{Tubular Neighborhood Theorem}
First we give a brief review of some related concepts from differential geometry.
For each $x\in \R^n$, the tangent space $T_x\R^n$ is canonically identified with $\R^n$
and the tangent bundle $T\R^n$ is canonically diffeomorphic to $\R^n\times \R^n$.
Let $M\subseteq \R^n$ be an embedded $k$-dimensional submanifold.
For each $x\in M$, the normal space to $M$ at $x$ is defined to be the $(n-k)$-dimensional subspace $N_xM \subseteq T_x\R^n$
consisting of all vectors orthogonal to $T_xM$ with respect to the Euclidean inner product.
The normal bundle $NM:=\{(x, v)\mid x\in M, v\in N_xM\}$ of $M$ is an embedded $n$-dimensional submanifold of $T\R^n$ by
Theorem 6.23 in \cite{Lee2013}.

We define a smooth map $E: NM \rightarrow \R^n$ by $E(x,v) = x+v$ where $x\in M$ and $v\in N_xM$.
A \textbf{tubular neighborhood} of $M$ is a neighborhood $U$ of $M$ in $\R^n$ which is the diffeomorphic image under $E$ of an open subset $V\subseteq NM$
of the form $V = \{(x,v)\in NM:   \; \|v\|_2< \delta(x)  \}$ for some positive continuous function $\delta: M\rightarrow \R$.

\begin{theorem}[Tubular Neighborhood Theorem \cite{Lee2013}]
  \label{thm:tubular}
  Every embedded submanifold of $\R^n$ has a tubular neighborhood.
\end{theorem}

Tubular neighborhood theorem is also true for complex analytic manifold (See for instance Theorem 6.2 in \cite{Zeng2020}).

\begin{lemma}\label{lem:regular2}
Let $h= \{h_1,...,h_k\}\in\R[x]$ with $V_{\R}(h)\neq\emptyset$, where $k<n$.
Suppose that the Jacobian matrix $\Jac$ of $h$ attains full rank at any point of $V_{\R}(h)$.
Then there exists a tubular neighborhood $U$
of $V_{\R}(h)$ such that every point $b\in U$ has a unique projection $x_b\in V_{\R}(h)$ of minimum distance to $b$. Moreover,
$x_b$ is the projection of a unique isolated simple zero of the Lagrangian system of
\begin{eqnarray}\label{eq:optd3}
 \min \; & ( \|x- b\|_2^2  ) /2 \\
 s.t. \; & h(x)=0 \nonumber.
\end{eqnarray}
\end{lemma}
\proof
The Lagrangian system of (\ref{eq:optd3}) is as below:
\begin{eqnarray}\label{eq:normal}
& \left(
  \begin{array}{c}
    x_1-b_1 \\
    \vdots \\
    x_n-b_n \\
  \end{array}
\right) = \left(
            \begin{array}{ccc}
              \partial h_1/\partial x_1 & \cdots & \partial h_k/\partial x_1 \\
              \vdots & \ddots & \vdots \\
              \partial h_1/\partial x_n & \cdots & \partial h_k/\partial x_n \\ 
            \end{array}
          \right)_{n\times k }\cdot \; \left(
                         \begin{array}{c}
                           \lambda_1 \\
                           \vdots \\
                           \lambda_k \\
                         \end{array}
                       \right)\\
        & h(x) = 0. \nonumber
\end{eqnarray}


Since the Jacobian matrix of $h$ is of full rank at any point of $V_{\R}(h)$, the set $V_{\R}(h)$ is a  submanifold of $\R^n$.
Theorem~\ref{thm:tubular} guarantees that there is a tubular neighborhood of $V_{\R}(h)$ such that $b\in U$.
Next we show that $x_b$ is unique.
Suppose that there is another different point $x'_b\in V_{\R}(h)$ of minimum distance to $b$.
Let $v'= b-x'_b$ and $v= b-x_b$. 
Since $\Jac^t$ is of constant rank $k$, its column vectors form a basis of the normal space of dimension $k$.
So $v, v'$ are normal vectors of $x_b,x'_b$  respectively  and  $b=E(x_b,v) = E(x'_b,v')$. Since $b \in U$, $E$ is injective and it implies
$(x_b,v)=(x'_b,v')$. It is a contradiction. 

On the other hand, because the Jacobian matrix of $h$ is of full rank at any point of $V_{\R}(h)$,
the value of $\lambda=(\lambda_1,\ldots,\lambda_k)$ is uniquely defined by a zero $x$ of $h(x)$.
So $x_b$ is the projection of a unique zero $(x_b,\lambda_b)$ of Eqs.~(\ref{eq:normal}).

Next we show that $(x_b,\lambda_b)$ is an isolated simple zero of Eqs.~(\ref{eq:normal}).
Otherwise, there must exist another point $x'_b\in V_{\C}(h)$ near $x_b$ also satisfying (\ref{eq:normal}),
which is impossible by applying the Tubular Neighborhood Theorem for the complex analytic manifold $V_{\C}(h)$ following a similar argument as
we did for the real case.
\foorp

\begin{cor}\label{cor:regular3}
Let $f= \{f_1,...,f_k\}$ be a set of polynomials in the ring $\R[x]$ and $h=\{f_1+w_1/\sqrt{\beta},..., f_k+w_k/\sqrt{\beta}\}$ where $\beta$ is a positive constant.
Then $V_{\R}(h)$ has a tubular neighborhood $U \subset \R^{n+k}$. In addition, for any point $(\randpoint,b) \in U$,
there is a unique point $(x^*,w^*) \in V_{\R}(h)$ to minimize the distance from $(\randpoint,b)$ to $V_{\R}(h)$ and
$(x^*,w^*)$ is the projection of a unique isolated simple zero of 
the corresponding Lagrangian system of
\begin{eqnarray}\label{eq:optd2}
 \min \; & ( \|w-b\|_2^2 + \|x- \randpoint\|_2^2  ) /2 \\
 &s.t. \hspace{1cm} h=0 \nonumber.
\end{eqnarray}
\end{cor}
\proof
Note that  $h$ attains full rank at any point of $V_{\R}(h)$.
Then the conclusion follows directly from Lemma~\ref{lem:regular2}.
\foorp




\begin{theorem}
Let $f= \{f_1,...,f_k\}$ be a set of polynomials in the ring $\R[x_1,...,x_n]$ with  $V_{\R}(f)\neq \emptyset$.
For any $\beta>0$, there exists an open set in $\R^n$ containing $V_{\R}(f)$.
And for any point $\randpoint$ in this set there is a unique point on $V_{\R}(h)$ to minimize the distance from $(\randpoint,0)\in\R^{n+k}$ to $V_{\R}(h)$,
where $h=\{f_1+w_1/\sqrt{\beta},..., f_k+w_k/\sqrt{\beta}\}$.

Moreover, if $\randpoint$ moves along some piecewise smooth curve $\mathcal{C}$ in this set, the corresponding projection trajectory can be obtained by solving Equation (\ref{eq:newsys3}) continuously from an initial point $\randpoint_0 \in \mathcal{C}$.
\end{theorem}
\proof
By Corollary~\ref{cor:regular3}, $V_{\R}(h)$ has a tubular neighborhood $U$ containing $V_{\R}(h)$ in $\R^{n+k}$
such that for any point $(\randpoint,b) \in U$,
there is a unique point $(x^*,w^*) \in V_{\R}(h)$ to minimize the distance from $(\randpoint,b)$ to $V_{\R}(h)$ and
$(x^*,w^*)$ is the projection of a unique isolated simple zero of 
the corresponding Lagrangian system of Equation~(\ref{eq:optd2}):
\begin{equation}\label{eq:normal2}
\left(
  \begin{array}{c}
    x_1-\randpoint_1 \\
    \vdots \\
    x_n-\randpoint_n \\
    w_1-b_1\\
     \vdots \\
    w_k-b_k \\
  \end{array}
\right) = \left(
            \begin{array}{ccc}
              \partial f_1/\partial x_1 & \cdots & \partial f_k/\partial x_1 \\
              \vdots & \ddots & \vdots \\
              \partial f_1/\partial x_n & \cdots & \partial f_k/\partial x_n \\
              1/\sqrt{\beta} &  &   \\
                & \ddots &  \\
                &        & 1/\sqrt{\beta} \\
            \end{array}
          \right)_{(n+k)\times k }\cdot \; \left(
                         \begin{array}{c}
                           \lambda_1 \\
                           \vdots \\
                           \lambda_k \\
                         \end{array}
                       \right).
\end{equation}
with $h=0$.
Since $\emptyset \neq V_{\R}(f) = \pi_{x}(V_{\R}(h) \cap V_{\R}(w)) \subset  \pi_{x}(U \cap   V_{\R}(w))$
and $U$ is an open set,
letting $U' = \pi_{x}(U \cap   V_{\R}(w))$,
we know that $U'$ is an open set in $\R^n$ containing $V_{\R}(f)$.
For any point $\randpoint \in U'$, we have  $(\randpoint, 0)\subset U$.
Thus, by Corollary~\ref{cor:regular3} there is a unique point $(x^*,w^*) \in V_{\R}(h)$ to minimize the distance from $(\randpoint,0)$ to $V_{\R}(h)$ and
$(x^*,w^*)$ is the projection of a unique isolated simple zero of  Equation~(\ref{eq:normal2}).
Since there is a one-to-one correspondence between the zeros of Equation~(\ref{eq:normal2})
and those of Equation (\ref{eq:newsys3}), $x^*$ is an isolated simple zero of Equation (\ref{eq:newsys3}).
Therefore, the projection trajectory can be obtained by solving Equation (\ref{eq:newsys3}) continuously from an initial point $\randpoint_0 \in \mathcal{C}$. 
\foorp


 A natural question is how to move  $\randpoint$ into the tubular neighborhood.
Recall the equation
\begin{equation}\label{eq:square}
x + \beta\cdot \Jac^t \cdot f= \randpoint.
\end{equation}
Suppose $x_0$, $\randpoint_0$ and $\beta_0$ satisfy this equation.
So $\lambda (x_0-\randpoint_0) +  \lambda \beta_0\cdot \Jac^t \cdot f= 0$ for any $\lambda>0$.
Let $\randpoint_1 = (1-\lambda)x_0+ \lambda \randpoint_0$  and $\beta_1= \lambda \beta_0$.
Then $(x_0-\randpoint_1)+ \beta_1\cdot \Jac^t \cdot f= 0$.

So $x_0$ is a critical point of Equation (\ref{eq:square}) with $\randpoint= \randpoint_1$ and $\beta=\beta_1$.
When $\lambda$ is small, $\randpoint_1$ is very close to $x_0$. By Theorem \ref{thm:est2}, $x_0$ should be quite close to the curve for sufficiently large $\beta$  and consequently $\randpoint_1$ is also close to the curve.
Since $\beta_1= \lambda \beta_0$, we can increase the value of $\beta$ by the homotopy (\ref{eq:homotopy2}) and it gives a refined solution
$x_1$ even closer to the curve.

In summary, we have the following algorithm.

\begin{algorithm}
\LinesNumbered
\caption{${\sf Move\;towards\;Tubular\;Neighborhood}$}
\label{Algo:MTN}
\KwIn{
  A system $f =\{f_1,\ldots,f_{k}\}\in\R[x_1,\ldots,x_n]$.
  A contraction factor $0<\lambda<1$. An initial point $(x_0,\randpoint_0)$ satisfying Equation (\ref{eq:square}) with a fixed penalty factor $\beta_0 \gg 1$.
}
\KwOut{
  A point $(x_1,\randpoint_1)$ satisfying Equation (\ref{eq:square}), where $\beta=\beta_0$, with $\randpoint_1 = (1-\lambda)x_0+ \lambda \randpoint_0$.
}
\Begin{
    let $\randpoint_1 = (1-\lambda)x_0+ \lambda \randpoint_0$ and $\beta_1= \lambda \beta_0$\;
    
    let $t_1= 1/ \beta_1$ and $t_0 = 1/ \beta_0$\;
    
    construct $H(x,t) = t(x-\randpoint_1)+ \Jac^t\cdot f \equiv 0$ with an initial point $(x_0,t_1)$\;
    
    move the parameter $t$ from $t_1$ to $t_0$ continuously yields $H(x_1,t_0)=0$\;
    
   return $(x_1,\randpoint_1)$. 
}
\end{algorithm}

\begin{remark}
Note that Algorithm \ref{Algo:MTN} cannot guarantee that $\randpoint_1$ belongs to the tubular neighborhood of $V_{\R}(f)$
since there is no information about the radius of the tubular neighborhood and the contraction factor is just chosen by the user.

However, the tubular neighborhood condition is sufficient but unnecessary for regular Jacobian of Equation (\ref{eq:square}).
Therefore, we can monitor the Jacobian during changing $\randpoint$ along some direction. If the Jacobian is close to being singular, it indicates  
that
$\randpoint$ might be out of the tubular neighborhood and we can use Algorithm \ref{Algo:MTN} to move $\randpoint$
towards the tubular neighborhood. Since the tubular neighborhood is an open set, $\randpoint$  will be in the neighborhood after finitely many iterations. Thus, this
restores the regularity of Equation (\ref{eq:square}).
\end{remark}

\subsection{Tracing direction}

The next question to answer is how to choose the direction to move $\randpoint$ in Equation (\ref{eq:square}).
By our assumption, $V_{\R}(f)$ is a one dimensional curve and the best choice of moving $\randpoint$ is along the tangent direction
of the curve. However it will be difficult to obtain such a direction if Jacobian of $f$ is rank-deficient along the whole curve.
Here we present a method for finding the tracing direction by computing eigenvalues.

Theorem \ref{thm:tubular} gives the geometric relation between $\randpoint$ and $x$.
That is there exists a $w$ such that $(x, w)$ is the projection of $(\randpoint,0)$
onto $V_{\R}(h)$, where $h = f +w/\sqrt{\beta}$.
When $\beta$ is sufficiently large, by Theorem \ref{thm:est2} we can consider $x$ as the projection of $\randpoint$ onto $V_{\R}(f)$ approximately.
If $\randpoint$ is quite close to the curve $V_{\R}(f)$, roughly speaking when the change of $\randpoint$ is
parallel to the tangent at $x$, then  we should have $\Delta \randpoint \approx  c \Delta  x$ for some constant $c$ close to $1$.
This implies that
\begin{equation}
(\beta\frac{\partial \Jac^t \cdot f}{\partial x} + I) \Delta  x = \Delta \randpoint \approx  c \Delta  x.
\end{equation}
Therefore, $\Delta \randpoint$ will be approximately equal to the eigenvector corresponding to the eigenvalue $1$.

Since $V_{\R}(f)$ is of dimension one and this eigenvector is close to the tangent, other eigenvectors will approximately span   the  normal space
of the curve at this point. Intuitively, a large change of  $\randpoint$ in the approximate normal space will not change the projection
$x$ too much and consequently the eigenvalues of such eigenvectors will be very large.

\begin{prop}
The matrix $\beta\frac{\partial \Jac^t \cdot f}{\partial x} + I$ is symmetric and  it has $n$ orthogonal eigenvectors.
\end{prop}
\proof
We only need to show the matrix  $A = \frac{\partial \Jac^t \cdot f}{\partial x}$ is symmetric.
Since $f = (f_1,...,f_k)$, it is straightforward to verify that
\begin{equation}
A_{ij} = \frac{\partial (\sum_{\ell=1}^k  f_{\ell} \frac{\partial f_{\ell}}{\partial x_i} )  }{\partial x_j}=
\sum_{\ell=1}^k (\frac{\partial f_{\ell}}{\partial x_i}\cdot \frac{\partial f_{\ell}}{\partial x_j}+ f_{\ell}  \frac{\partial^2 f_{\ell}}{\partial x_i\, \partial x_j}  ) = A_{ji}.
\end{equation}

The eigenvalues of a real symmetric matrix are real and their eigenvectors are orthogonal.
\foorp \\

{
\begin{example}
  Recall the example $f:=(x_1^3-x_2)^2$ in Section~\ref{sec:intro-ex}.
  With $\beta=10^4$, at the witness point $(-0.8296,-0.5982)$,
  we have
  $$
  \beta\frac{\partial \Jac^t \cdot f}{\partial x} + I=\left(\begin {array}{cc}  188.7722&- 91.9182
\\ \noalign{\medskip}- 91.9182& 45.5187\end {array}
  \right),
  $$
  which has two eigenvalues $233.6758$ and $0.6150$,
  with corresponding eigenvectors $(0.8985, -0.4389)$ and $(0.4389, 0.8985)$.
  So the second eigenvector, whose corresponding eigenvalue close to $1$, will be chosen as the tracing direction,
  which is indeed an approximation of the tangent direction $(0.4359, 0.9000)$ of the curve at the witness point.
\end{example}
}

\subsection{Companion curve tracing method}
We are now ready to present the companion curve tracing method in detail.
The input is a finite set $f$ of polynomials in $\R[x]$.
The algorithm starts by applying the criterion in Section~\ref{s:emptiness}
to determine if $V_{\R}(f)$ is empty. If $V_{\R}(f)$  is determined to empty,
the algorithm terminates and return $\emptyset$.
Otherwise, it will produce given number of points in $V_{R}(h)$, where $h=\{f_1+w_1/\sqrt{\beta}, \ldots, f_k+w_k/\sqrt{\beta}\}$
for a large $\beta\gg 1$, such that $\|f\|\leq \epsilon$ holds.
More precisely, we have the following algorithm.
\begin{itemize}
\item Algorithm {\sf CompanionCurveTracing}
\item Input:
  \begin{itemize}
  \item a finite set of polynomials $f=\{f_1,\ldots,f_k\}$.
  \item a prescribed number $N$.
  \end{itemize}
\item Output: return $\emptyset$ if $V_{\R}(f)$ is determined to be empty by the algorithm;
  otherwise return at least $N$ points for each connected component of $V_{\R}(f)$ such that for each point $p$, the backward error $\|f(p)\|\leq\epsilon$.
\item Steps:
\begin{enumerate}
\item choose a random point $\randpoint$ in $\R^n$ and a large number $\beta \gg 1$
\item apply Proposition~\ref{prop:empty} to determine if $V_{\R}(f)$ is empty; if true, then return $\emptyset$\;
    \item obtain the real solution set $S$ of the square system $(x-\randpoint)+ \beta \Jac^t\cdot f =0$ by homotopy continuation method
    \item let $S' = \{x:   \|f(x)\| < \epsilon  \}$ for some tolerance $\epsilon$
    \item if $S'=\emptyset$, then return $\emptyset$
    \item for each $x \in S'$
      \begin{enumerate}
        \item[6.1] call Algorithm \ref{Algo:MTN} to get a new point $(x',\randpoint')$ where $(\randpoint', 0)$ is presumably in
           the tubular neighborhood of $V_{\R}(f+w/\sqrt{\beta})$ 
        \item[6.2] find the unit eigenvector $\Delta x$ of $\beta\frac{\partial \Jac^t \cdot f}{\partial x} + I$ at $x'$ with eigenvalue $c$ close to $1$
        \item[6.3] update $\randpoint' = \randpoint'+ h c\Delta x$ and $x' = x' + h \Delta x$, where $h$ is the step size
        \item[6.4] refine $x'$ by Newton iteration with fixed  $\randpoint'$ { by solving $(x'-\randpoint')+ \beta \Jac^t\cdot f =0$.}
        \item[6.5] if the smallest eigenvalue of $\beta\frac{\partial \Jac^t \cdot f}{\partial x} + I$ at $x'$ is close to zero, it indicates that $(\randpoint',0)$ may be out of the tubular neighborhood, goto step 6.1
        \item[6.6] goto step 6.2 with updated $(x', \randpoint')$ until we have $N$ points of the curve starting from $x$
    \end{enumerate}
    \item goto step 6 until we enumerate all points in $S'$        
\end{enumerate}
\end{itemize}

\section{Examples}
\label{sec:ex}
In this section, we illustrate the effectiveness of the companion curve method on tracing
some curves defined by rank-deficient systems.
{ For all the examples below, we simply choose $\beta=10^4$.}

\begin{example}[\cite{Lax05}]
  \label{ex:lax}
  Let $f$ be the discriminant of the characteristic polynomial of the following matrix
  $$
  \left(
  \begin{array}{rcl}
    x_1 & 1 & 1\\
    1 & x_2 & 1\\
    1 & 1 & x_3\\
  \end{array}
  \right).
  $$
  It was proved that $f$ is a sum of squares~\cite{Lax05}, which is necessarily rank-deficient at any point of $V_{\R}(f)$.
  The companion curve method generates a straight line as illustrated by Fig.~\ref{fig:lax}.

   \begin{figure}
    \centering
 \begin{minipage}[t]{0.48\linewidth}
\centering
\includegraphics[width=\textwidth]{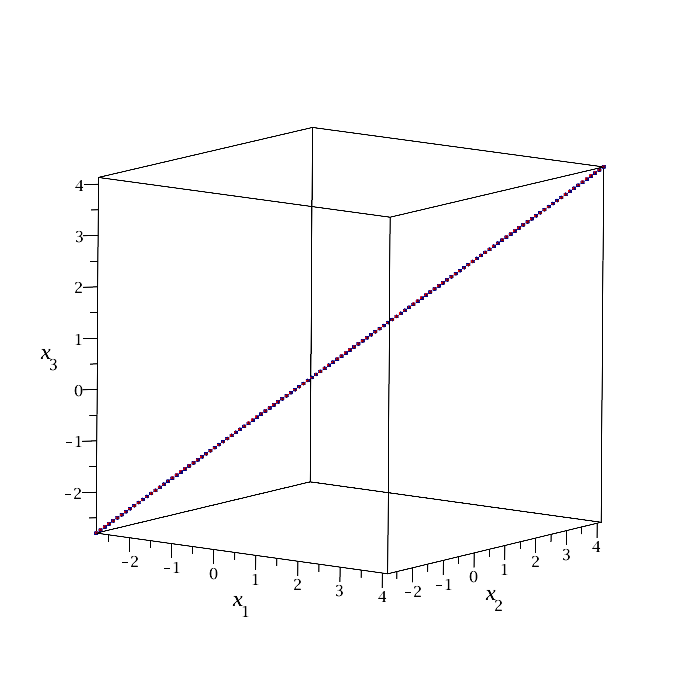}
\end{minipage}
\begin{minipage}[t]{0.48\linewidth}
\centering
\includegraphics[width=\textwidth]{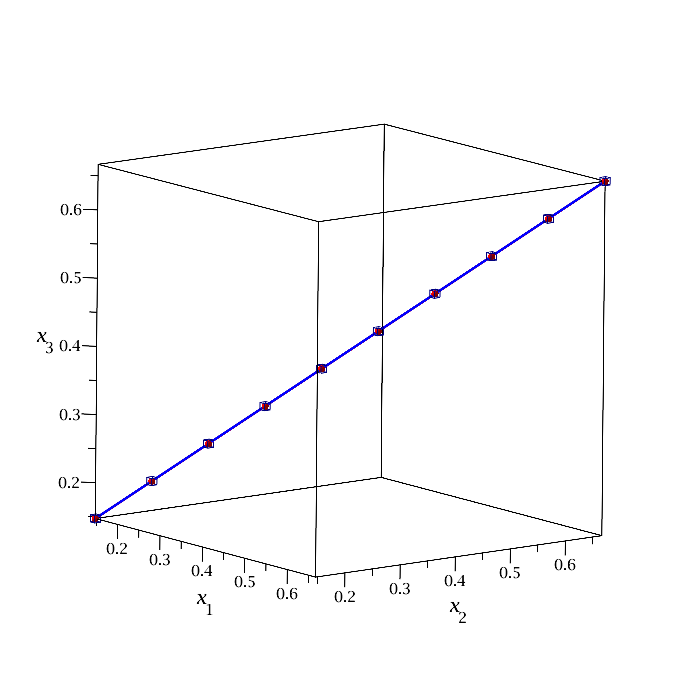}
\end{minipage}
\caption{{(Color online) Left: The traced curve (in red $\bullet$) of $V_{\R}(f)$ in Example~\ref{ex:lax} and its companion curve (in blue $\Box$) overlap.
    Right: zooming in on the left figure to see the difference.}}
\label{fig:lax}
\end{figure} 
  
\end{example}

\begin{example}[Example $7$ in~\cite{Chen2017}]
\label{ex7}
Consider the following polynomial system
{\small
$$
f := \left\{
\begin{array}{rcl}
3\,{ x_4}\,{{ x_1}}^{3}&+&100\,{{ x_2}}^{4}-119\,{{ x_3}}^{2}{{x_2}}^{2}-6\,{{ x_4}}^{2}{ x_3}\,{ x_2}+36\,{{ x_3}}^{4}+9\,{{ x_4}}^{4}+3\,{ x_4}\,{{ x_1}}^{2}\\
                     &-&40\,{{ x_2}}^{2}{ x_1}-16\,{ x_3}\,{ x_2}\,{ x_1}+24\,{{ x_3}}^{2}{ x_1}+48\,{{x_4}}^{2}{ x_1}+71\,{{ x_1}}^{2}\\
                     &-&7\,{ x_2}\,{ x_1}-20\,{{x_2}}^{2}+2\,{ x_3}\,{ x_2}+12\,{{ x_3}}^{2}-6\,{{ x_4}}^{2}-8\,{ x_1}-7\,{ x_2}+3,\\
9\,{{ x_4}}^{4}&-&6\,{{ x_4}}^{2}{ x_3}\,{ x_2}+3\,{ x_4}\,{{ x_1}}^{3}+36\,{{ x_3}}^{4}-119\,{{x_3}}^{2}{{ x_2}}^{2}+100\,{{ x_2}}^{4}+48\,{{ x_4}}^{2}{ x_1}\\
              &-&3\,{ x_4}\,{{ x_1}}^{2}+24\,{{ x_3}}^{2}{ x_1}-16\,{ x_3}\,{x_2}\,{ x_1}-40\,{{ x_2}}^{2}{ x_1}-6\,{{ x_4}}^{2}+12\,{{x_3}}^{2}\\
              &+&2\,{ x_3}\,{ x_2}-20\,{{ x_2}}^{2}-7\,{ x_2}\,{ x_1}+71\,{{ x_1}}^{2}+7\,{ x_2}-14\,{ x_1}+1\\
\end{array}
\right\},
$$
}
which defines a real algebraic curve according to~\cite{Chen2017}. 
However, on the other hand, this system has two polynomials with four variables $x_1,x_2,x_3,x_4$.
{ Hence, the nullity of $\Jac_f$ is at least $2$,
which implies that  $f$ must be rank-deficient.
}
The companion curve method generates several curves, as illustrated by Fig.~\ref{fig:ex7}.

\begin{figure}
    \centering
 \begin{minipage}[t]{0.49\linewidth}
\centering
\includegraphics[width=\textwidth]{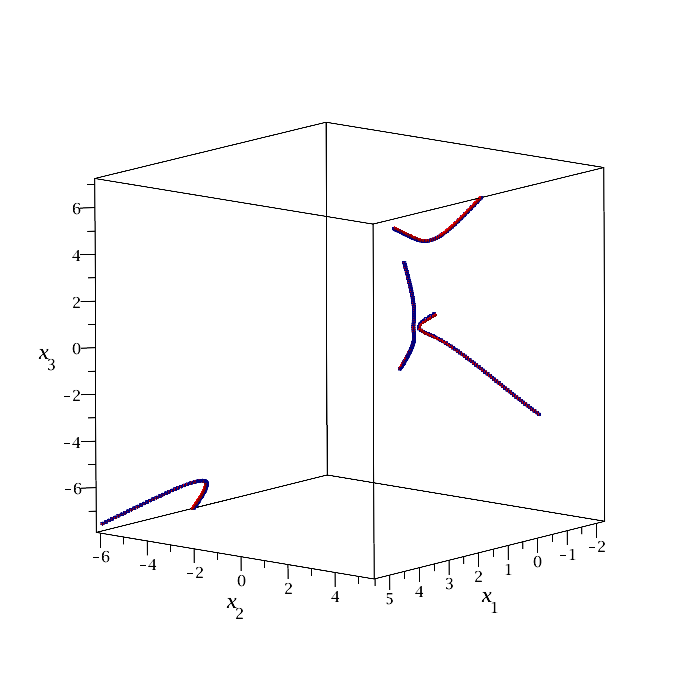}
\end{minipage}
\begin{minipage}[t]{0.48\linewidth}
\centering
\includegraphics[width=\textwidth]{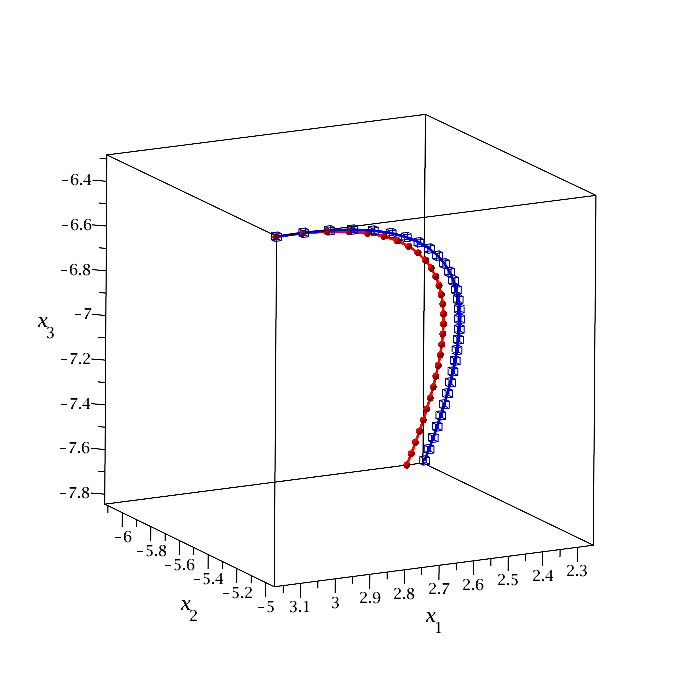}
\end{minipage}
\caption{{(Color online) Left: Projections in $(x_1, x_2, x_3)$ of the traced curve (in red $\bullet$) of $V_{\R}(f)$ in Example~\ref{ex7} and its companion curve (in blue $\Box$) overlap.
    Right: zooming in part of the left bottom branch of $V_{\R}(f)$ (and its companion curve) to see the difference.}}
\label{fig:ex7}
\end{figure} 

\end{example}

\begin{example}[Example 8 in~\cite{Chen2017}]
\label{ex8}
Consider the following system $f$:
{\small
$$
\left\{
\begin{array}{rcl}
{{ x_1}}^{3}&-&16\,{ x_3}\,{{ x_1}}^{2}+100\,{{ x_2}}^{2}{ x_1}+160\,{ x_4}\,{ x_2}\,{ x_1}+64\,{{ x_3}}^{2}{ x_1}+64\,{{ x_4}}^{2}{ x_1}-14\,{{ x_1}}^{2}\\
          &+&20\,{ x_2}\,{ x_1}+112\,{x_3}\,{ x_1}+16\,{ x_4}\,{ x_1}-4\,{ x_3}\,{ x_5}-3\,{x_4}\,{ x_5}+50\,{ x_1}-8\,{ x_3}+1,\\
{ x_2}\,{{ x_1}}^{3}&-&16\,{ x_3}\,{ x_2}\,{{ x_1}}^{2}+100\,{{ x_2}}^{3}{ x_1}+160\,{ x_4}\,{{ x_2}}^{2}{ x_1}+64\,{{ x_3}}^{2}{ x_2}\,{ x_1}\\
                 &+&64\,{{ x_4}}^{2}{ x_2}\,{ x_1}-14\,{ x_2}\,{{ x_1}}^{2}+20\,{{ x_2}}^{2}{ x_1}+112\,{ x_3}\,{ x_2}\,{ x_1}+16\,{ x_4}\,{ x_2}\,{ x_1}\\
                 &-&4\,{ x_2}\,{ x_3}\,{ x_5}-3\,{ x_2}\,{x_4}\,{ x_5}-{{ x_1}}^{2}+50\,{ x_2}\,{ x_1}+16\,{ x_3}\,{ x_1}-100\,{{ x_2}}^{2}-8\,{ x_3}\,{ x_2}\\
                 &-&160\,{ x_4}\,{x_2}-64\,{{ x_3}}^{2}-64\,{{ x_4}}^{2}+14\,{ x_1}-19\,{ x_2}-112\,{ x_3}-16\,{ x_4}-50,\\
-6\,{ x_4}\,{ x_1}&+&2\,{ x_2}-7\,{ x_5}+1,\\
-2\,{ x_6}\,{ x_1}&-&4\,{ x_4}\,{ x_3}+3\,{ x_5}+1
\end{array}
\right\},
$$
}
which also defines a real algebraic curve according to~\cite{Chen2017}. 
This system has four polynomials with six variables $x_1,\ldots, x_6$.
{ Hence, the nullity of $\Jac_f$ is at least $2$,
which implies that  $f$ must be rank-deficient.}
The companion curve method generates several curves, as illustrated by Fig.~\ref{fig:ex8}.

\begin{figure}
    \centering
 \begin{minipage}[t]{0.48\linewidth}
\centering
\includegraphics[width=\textwidth]{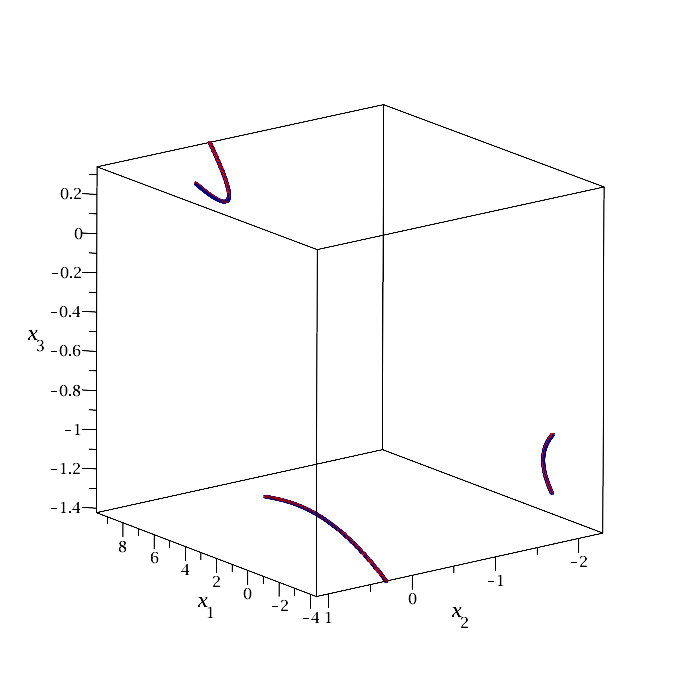}
\end{minipage}
\begin{minipage}[t]{0.48\linewidth}
\centering
\includegraphics[width=\textwidth]{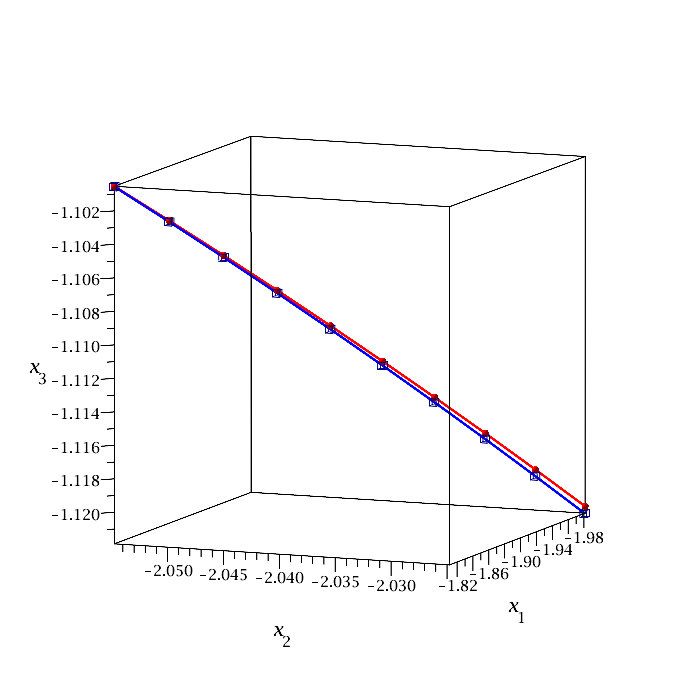}
\end{minipage}
\caption{{(Color online) Left: Projections in $(x_1, x_2, x_3)$ of the traced curve (in red $\bullet$) of $V_{\R}(f)$ in Example~\ref{ex8} and its companion curve (in blue $\Box$) overlap.
    Right: zooming in on part of the right bottom branch of $V_{\R}(f)$ (and its companion curve) to see the difference.}}
\label{fig:ex8}
\end{figure} 

\end{example}

\begin{example}
  The Choi-Lam polynomial \cite{Choi1977} $$f= {x}^{2}{y}^{2}+{x}^{2}{z}^{2}+{y}^{2}{z}^{2}-4\,xyz+1$$ has $4$ real isolated solutions: $\{(1,1,1),(1,-1,-1),(-1,1,-1),(-1,-1,1)\}$
  and is rank-deficient at all these points.
  The real zero set of $f$ can be seen as a degenerate curve.
  Fig.~\ref{fig:lam} illustrates the points and their corresponding companion curve produced by the method.
  Note that the companion curve bounces back and forth around the isolated solutions.
  This is due to the fact that in Algorithm {\sf CompanionCurveTracing} Step $6.3$ moves the guiding point $\randpoint$ away from
  an isolated point while Step $6.1$ will drag $\randpoint$ back to a neighborhood of the isolated point. 
  \begin{figure}
    \centering
 \begin{minipage}[t]{0.48\linewidth}
\centering
\includegraphics[width=\textwidth]{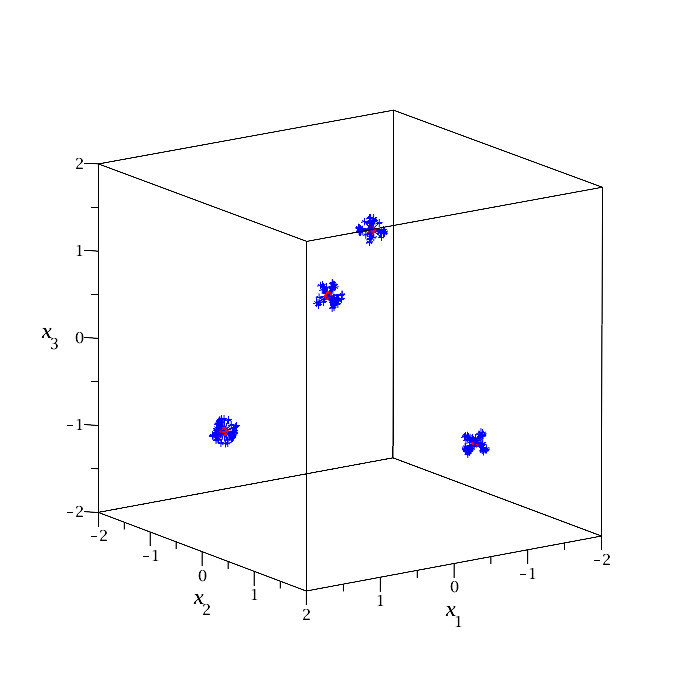}
\end{minipage}
\begin{minipage}[t]{0.48\linewidth}
\centering
\includegraphics[width=\textwidth]{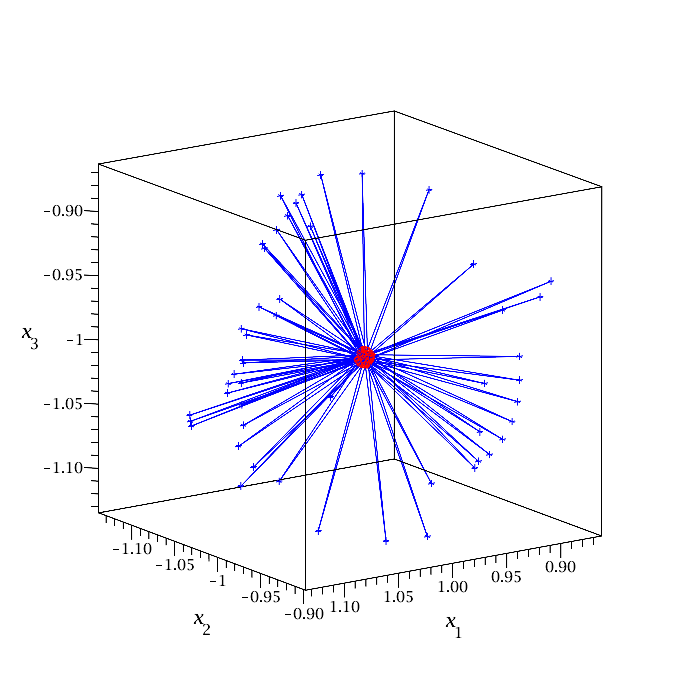}
\end{minipage}
\caption{{(Color online) Left: The four points traced (in red $\bullet$) of $V_{\R}(f)$ and their companion curve (in blue $+$).
    Right: zooming in on the left figure around point $(1, 1, 1)$.}}
\label{fig:lam}
\end{figure} 
\end{example}

\section{Conclusion and future work}
\label{sec:con}
In this paper, we proposed a companion curve method for tracing
real algebraic curves whose defining system has singular Jacobian
along the whole curve.
The effectiveness of this method is illustrated by several non-trivial
systems. We notice that, even the algebraic form of the real algebraic curve
is singular, as long as the curve is geometrically regular or has only singletons, the method works
very well. We do notice that the method may not be able to obtain a complete
tracing of  curves which are not only ``algebraically singular'' but also ``geometrically singular'', simple as $(xy)^2$,
without an algebraic preprocess.
This is due to the fact that a slight perturbation of the curve may generate more than one connected components although
the curve itself is connected via the singular point.
A preliminary investigation suggests that  the symmetry matrix $\frac{\partial \Jac^t \cdot f}{\partial x} + I$
plays a vital role in suggesting tracing directions,
which deserves a further in-depth investigation.

\medskip\noindent{\bf \normalsize Acknowledgements.}
This work was funded by NSFC (No. 11771421), the Key Research Program of Frontier Sciences of CAS (No. QYZDB-SSW-SYS026), CAS ``Light of West China'' Program, National Key Research and Development Program (No. 2020YFA0712300), and Chongqing Programs (No. cstc2018jcyj-yszxX0002, No. cstc2019yszx-jcyjX0003, No. cstc2020yszx-jcyjX0005).

\end{document}